\documentclass[graybox]{svmult}

\usepackage{type1cm}        
\usepackage{makeidx}         
\usepackage{graphicx}        
\usepackage{multicol}        
\usepackage[bottom]{footmisc}

\usepackage{newtxtext}       %
\usepackage[varvw]{newtxmath}       

\usepackage{latexsym}
\usepackage{graphicx}
\usepackage{color}
\usepackage{bm}
\usepackage{algorithm, algorithmic}

\def\ub {{\bm u}}

\def\0b {{\bf 0}}

\def\ub {{\bf u}}

\def\ub {\bm{u}}

\def\nb {\bm{n}}

\def\ub {\bm{u}}

\def\nb {\bm{n}}

\makeindex             

\begin{document}

\title*{Biot model with generalized eigenvalue problems for scalability and robustness to parameters}
\titlerunning{Biot model with generalized eigenvalue problems} 
\author{Pilhwa Lee}
\institute{Department of Mathematics, Morgan State University, 1700 E. Cold Spring Lane, Baltimore, MD, USA ~~~\email{Pilhwa.Lee@morgan.edu}}
%
%
\maketitle

\abstract{We consider Biot model with block preconditioners and generalized eigenvalue problems for scalability and robustness to parameters. A discontinuous Galerkin discretization is employed with the displacement and Darcy flow flux discretized as piecewise continuous in $P_1$ elements, and the pore pressure as piecewise constant in the $P_0$ element with a stabilizing term. Parallel algorithms are designed to solve the resulting linear system. Specifically, the GMRES method is employed as the outer iteration algorithm and block-triangular preconditioners are designed to accelerate the convergence. In the preconditioners, the elliptic operators are further approximated by using incomplete Cholesky factorization or two-level additive overlapping Schwartz method where coarse grids are constructed by generalized eigenvalue problems in the overlaps (GenEO). Extensive numerical experiments show a scalability and parametric robustness of the resulting parallel algorithms.}

\section{Introduction}
\indent \indent Poroelasticity, {\it i.e.}, elasticity of porous media with permeated Darcy flow, is pioneered by Biot \cite{Biot1941general, Biot1955theory}. In this paper, we propose a numerical scheme for solving the Biot model with three-fields linear poroelasticity. We consider a discontinuous Galerkin discretization, {\it i.e.}, the displacement and Darcy flow flux discretized as piecewise continuous in $P_1$ elements, and the pore pressure as piecewise constant in the $P_0$ element with a stabilizing term. The emerging formulation is a saddle-point problem, and more specifically, a twofold saddle-point problem.\\ 
\indent  This indefinite system is computational challenging with slow convergence in iterative methods. It is necessary to incorporate relevant preconditioners for saddle-point problems. There have been decomposition methods through overlapping Schwarz methods \cite{Klawonn1998overlapping, Pavarino2000indefinite, Cai2015overlapping, Heinlein2016}.
~~We use GMRES as the outer iterative solver accelerated by a parallelized block-triangular preconditioner with overlapping additive Schwarz method (OAS) for displacement, and Darcy flow flux and Schur complement for pressure by Cholesky factorizations. In order to make the scheme scalable and robust to broad ranges of parameters and their potential heterogeneous distributions, the coarse grid should be well constructed. In the paper, we take the approach of constructing coarse spaces with eigenfunctions based on generalized eigenvalue problems \cite{Bjorstad2001, Galvis2010, Spillane2011}. Specifically, we devise a parallel preconditioner of two-level OAS with coarse grid construction by Generalized Eigenvalue problems in the Overlaps (GenEO) \cite{Spillane2014}.
\section{Linear poroelastic model}
Poroelastic models describe the interaction of fluid flows and deformable elastic porous media saturated in the fluid. Let $\bm{u}$ be the elastic displacement, $p$ be the pore pressure.
We assume that the permeability is homogeneous: ${\bf K} = \kappa{\bf I}$. Denote $\bm{z}$ as the Darcy volumetric fluid flux. The quasi-static Biot model reads as:
\begin{eqnarray}
-  \nabla \cdot \left(\sigma(\ub) - \alpha p {\bf I}  \right) &=& \bm{f},  \label{lee1}\\
{\bf K}^{-1} \bm{z} + \nabla p &=& \bm{b},   \label{lee2}\\
\frac{\partial}{\partial t} \left(\alpha \nabla \cdot \bm{u} + c_0 p \right)+ \nabla \cdot \bm{z} &=& g, \label{lee3}
\end{eqnarray}
where $\sigma(\ub)$ is the deviatoric stress, $\bm{f}$ is the body force on the solid, $\bm{b}$ is the body force on the fluid, $g$ is a source or sink term, $c_0>0$ is the constrained specific storage coefficient, $\alpha$ is the Biot-Willis constant which is close to 1. For the ease of presentation, we consider mixed partial Neumann and partial Dirichlet boundary conditions. Specifically, the boundary $\partial\Omega$ is divided into the following:
$$
\partial\Omega=\Gamma_{\rm d}\cup\Gamma_{\rm t} \quad ~\mbox{and}~ \quad \partial\Omega=\Gamma_{\rm p}\cup\Gamma_{\rm f},
$$
where $\Gamma_{\rm d}$ and $\Gamma_{\rm t}$ are for displacement and stress boundary conditions; $\Gamma_{\rm p}$ and $\Gamma_{\rm f}$ are for pressure and flux boundary conditions. ~
Accordingly, the boundary conditions are the following: 
\begin{eqnarray} \label{BC_mix}
\ub &=& \bm{0} \quad \mbox{on} ~\Gamma_{\rm d}, ~~~~~~~(\sigma(\ub)-\alpha p \mathbf{I}) \cdot \nb = \mathbf{t} \quad \mbox{on} ~\Gamma_{\rm t}, \\
p &=& 0\quad \mbox{on} ~\Gamma_{\rm p}, ~~~~~~~\bm{z} \cdot \nb = g_2 \quad \mbox{on} ~\Gamma_{\rm f}.
\end{eqnarray}
For simplicity, the Dirichlet conditions are assumed to be homogeneous.

%
%

\section{Saddle-point problem: discretization of $\mathbf{P}_1-\mathbf{P}_1-P_0$}
We apply the finite element method where domains are normally shaped as triangles in $\mathbb{R}^{2}$. Let $\mathcal{T}_{h}$ be a partition of $\Omega$ into non-overlapping elements $K$. We denote by $h$ the size of the largest element in $\mathcal{T}_{h}$.
On the given partition $\mathcal{T}_{h}$  we apply the following finite element spaces \cite{Berger2015}.
%
%
\begin{eqnarray}
\bm{V}_h  &:=& \{ \bm{u}_h \in (C^0(\Omega))^d: \mathbf{u}_h|K \in {\bf P}_1(K) ~\forall K \in \mathcal{T}^h, \bm{u}_h = 0~ {\rm on}~ \Gamma_{\rm d} \}  \\
\bm{W}_h  &:=& \{ \bm{z}_h \in (C^0(\Omega))^d: \mathbf{z}_h|K \in {\bf P}_1(K) ~\forall K \in \mathcal{T}^h, \bm{z}_h \cdot \bm{n} = 0~ {\rm on}~ \Gamma_{\rm f} \}  \\
Q_h &:=& \{ p_h : p_h|K \in {\bf P}_0(K) ~ \forall K \in \mathcal{T}^h \}
\end{eqnarray}
The problem is to find $(\bm{u}_h^n, \bm{z}_h^n, p_h^n) \in \bm{V}_h \times \bm{W}_h \times Q_h$ at the time step $n$~~such that
\begin{equation}
\left\{
\begin{array}{l}
a(\bm{u}_h^n, \bm{v}_h)  -(p_h^n, \nabla \cdot \bm{v}_h) = (\bm{f}^n, \bm{v}_h) + (\bm{t}^n,  \bm{v}_h)_{\Gamma_{\rm t}}, ~\forall \bm{v}_h \in \bm{V}_h   \\
(K^{-1} \bm{z}_h^n, \bm{w}_h) - (p_h^n, \nabla \cdot \bm{w}_h) = (\bm{b}^n, \bm{w}_h), ~\forall \bm{w}_h \in \bm{W}_h      \\
(\nabla \cdot \bm{u}_{\Delta t,h}^n, q_h) + (\nabla \cdot \bm{z}_h^n, q_h)+\frac{c_0}{\alpha}(p_h^n, q_h)  + J(p_{\Delta t,h}^n,q_h) = \frac{1}{\alpha} (g^n, q_h), ~\forall q_h \in Q_h
\end{array}
\right.
\end{equation}
where
\begin{displaymath}
J(p,q) = \delta_{\rm STAB} \sum_K \int_{\partial K \backslash \partial \Omega} h_{\partial K} [p][q]ds
\end{displaymath}
is a stabilizing term \cite{Burman2007}, and $p_{\Delta t,h}^n = (p_h^n - p^{n-1}_h)/\Delta t$. The finite element discretization will lead to a twofold saddle-point problem of the following form:


\begin{equation}
\left[\begin{array}{ccc}
A_{\bm{u}}  & 0      &  B^{T}_{1} \\[1mm]
0  & A_{\bm{z}}      & B^{T}_{2} \\[1mm]
B_1        &B_2    & -A_p \\
\end{array}\right]
\left[\begin{array}{c}
\bm{u}_h \\[1mm]
\bm{z}_h \\[1mm]
p_h
\end{array}\right]=
\left[\begin{array}{c}
{\bm f}_1 \\[1mm]
{\bm f}_2\\[1mm]
{\bm f}_3
\end{array}\right]. \label{matrix_A}
\end{equation}

%
\noindent Denote the block matrices of $A$, $B$, and the Schur complement $S$ in the following:
\begin{equation}
A=\left[\begin{array}{cc}
A_{\bm{u}}  & 0       \\
0  & A_{\bm{z}}
\end{array}
\right], \quad
B= \left[\begin{array}{cc}
B_1        &B_2
\end{array}
\right], \quad S=-(BA^{-1}B^T + A_p).
\end{equation}
Usually, for saddle-point problem of the form (\ref{matrix_A}), one takes the preconditioner as a block lower-triangular, 
\begin{equation}
T=\left[\begin{array}{cc}
A    & 0 \\[1mm]
B   & S      \\[1mm]
\end{array}\right].
\end{equation}

\subsection{Two-Level Additive Schwarz algorithm (OAS-2) for $A_{\bm{u}}$}

We now introduce the decomposition into local and coarse spaces.
The local problems are defined on the extended subdomains $\Omega_i^\prime$. To each of the $\Omega_i^\prime$, we associate a local space
\begin{equation}
\bm{V}_i =  \bm{V}^h(\Omega_i^\prime)
\cap \bm{H}^1_0(\Omega_i^\prime),
\end{equation}
and a bilinear form $a_i^\prime(\bm{u}_i,\bm{v}_i) := a(R_i^T\bm{u}_i,R_i^T \bm{v}_i)$, where $R_i^T: \bm{V}_i \rightarrow  \bm{V}^h,$ simply extends any element of $\bm{V}_i$ by zero outside $\Omega_i^\prime$. Then, as we will only consider algorithms for which the local problems are solved exactly, we find that the local operators are
\begin{equation}
A_i^{\prime} = R_i A R_i^T, \quad i=1, ..., N.
\end{equation}
Given the local and coarse embedding operators ${R}_{i}^T:\bm{V}_{i}\rightarrow \bm{V}^h$, $i=1,...,N$, and ${R}_0^T:\bm{V}_{0}\rightarrow \bm{V}^h$, the discrete space $\bm{V}^h$ can be decomposed into coarse and local spaces as
\begin{equation}
\bm{V}^h = {R}_0^T \bm{V}_0   + \sum_{i}{R}_{i}^T\bm{V}_{i}.
\end{equation}
%
%
The coarse space on the coarse subdomain mesh $\tau_H$ is denoted by
\begin{equation}
\bm{V}_0 = \bm{V}^H := \{\bm{v}\in \bm{V} : \bm{v}|_{\Omega_i} \in ({\bf P}_1(\Omega_i))^d~~ \forall \Omega_i \in \tau_H \}.
\end{equation}


\subsection{Construction of coarse spaces by GenEO for $A_u$}
For all subdomains $1 \leq i \leq N$, the generalized eigenvalue problem is to find $(V_{ik}, \lambda_{ik}) \in$ range$(A^\prime_i) \times \mathbb{R}$ such that
\begin{equation}
 A^{\prime}_i V_{ik} = \lambda_{ik} D_i A^{\prime}_i D_i V_{ik}.
\end{equation}
where $\{ D_i \}_{i=1}^N$ defines a partition of unity, $\sum_{i=1}^N R_i^T D_i R_i = I$.
The GenEO coarse space $V_0$ is based on the following local contributions:
\begin{equation}
 Z_{i\tau} = {\rm span}\{V_{ik} | \lambda_{ik} < \tau\},
\end{equation}
which are weighted with the partition of unity as follows:
\begin{equation}
V_0 = \oplus_{i=1}^N R_i^T D_i Z_{i \tau}.
\end{equation}
When $Z_0$ be a column matrix so that $V_0$ is spanned by its columns, $R_0 = Z_0^T$. Two-level overlapping Schwarz method by GenEO is summarized in Algorithm 1.
%

\begin{algorithm}[hbt]
\caption{Two-level overlapping Schwarz method by GenEO \cite{Jolivet2021}}
\begin{algorithmic}
\STATE{1. Solve the local generalized eigenvalue problem:}
 \STATE{~~~~~~$A^{\prime}_i y_i = \lambda_i \tilde{R}_i R_i^T  A^{\prime}_i R_i \tilde{R}_i^T y_i.$}
\STATE{~~~~~~~$\{ \tilde{R}_i \}_{i=1}^N$ are the same operators as $\{ R_i \}_{i=1}^N$ except that entries on the overlap are set to 0 \cite{Cai1999}.}
\STATE{2. Collect the $\nu_i$ smallest eigenpairs $\{ y_{i_j}, \lambda_{i_j} \}^{\nu_i}_{j=1}$.}
\STATE{3. Assemble a local deflation dense matrix:}
\STATE{~~~~~~~$W_i = [ D_i y_{i_1} \cdots D_i y_{i_{\nu_i}} ].$}
\STATE{~~~~~~~$\{ D_i = \tilde{R}_i  R_i^T \}_{i=1}^N$ defines the partition of unity.}
\STATE{4. Define a global deflation matrix:}
\STATE{~~~~~~~$P = [ R_1^T W_1 \cdots R_N^T W_N].$}
\STATE{5. Define a two-level preconditioner using the Galerkin product of $A$ and $P$:}
\STATE{~~~~~~~$M^{-1} = \sum_{i=1}^N \tilde{R}_i^T (R_i A R_i^T)^{-1} R_i,$}
\STATE{~~~~~~~$Q = P(P^T AP)^{-1}P^T$,}
\STATE{~~~~~~~$M^{-1}_{\rm additive} = Q + M^{-1} ~~{\rm or}~~ Q + M^{-1}(I - AQ)$.}
\end{algorithmic}
\label{Monte_Carlo_method} \label{algorithm1}
\end{algorithm}

\section{Numerical experiments}

A test problem is formulated with $\alpha = 1$, $c_0 = 0$, $\Omega = [0,1]^2$ and $t \in [0, 0.25]$:
\begin{eqnarray}
-(\lambda + \mu) \nabla (\nabla \cdot \bm{u}) - \mu \nabla^2 \bm{u} + \nabla p &=& 0,  \label{test1} \nonumber \\
 \mathbf{K}^{-1} \bm{z} + \nabla p &=& 0,   \label{test2}\\
\nabla \cdot (\bm{u}_t + \bm{z}) &=& g_1. \label{test3} \nonumber
\end{eqnarray}
%

The involving initial and boundary conditions are the following:
\begin{equation}
\left\{
\begin{array}{l}
\ub~=\bm{0}\quad \mbox{on} ~\partial \Omega = \Gamma_{\rm d}, \\
\bm{z} \cdot\nb~=g_2 \quad \mbox{on} ~\partial \Omega = \Gamma_{\rm f},\\
\bm{u}(\bm{x},0) = 0, \bm{x} \in \Omega,\\
p(\bm{x},0) = 0, \bm{x} \in \Omega.
\end{array}
\right.
\end{equation}

%
%
We consider the following analytic solution
\begin{eqnarray}
\bm{u}  &=& \frac{-1}{4 \pi (\lambda + 2 \mu)} \left[\begin{array}{c}
\cos(2 \pi x) \sin(2 \pi y) \sin(2 \pi t)    \\
\sin(2 \pi x) \cos(2 \pi y) \sin(2 \pi t)
\end{array}
\right],  \nonumber \\
\bm{z}  &=& -2 \pi k \left[\begin{array}{c}
\cos(2 \pi x) \sin(2 \pi y) \sin(2 \pi t)    \\
\sin(2 \pi x) \cos(2 \pi y) \sin(2 \pi t)
\end{array}
\right], \\
p &=& \sin(2 \pi x) \sin(2 \pi y) \sin(2 \pi t), \nonumber
\end{eqnarray}
and derive the compatible source term of $g_1$.
%

\subsection{Numerical implementation}
%
%

As the focus of this paper is to justify the effectiveness and the efficiency of the algorithm, we mainly study the performance of the parallel preconditioner discussed as above. In our implementation, we use a fiinite element library, libMesh. We apply triangular element with 3 nodes. The GMRES method and overlapping Schwarz preconditioners are based on PETSc. The initial guess is zero and the stopping criterion is set as a $10^{-8}$, reduction of the residual norm. In each test, we count the iterations. For unstructured domain partition, we apply ParMETIS, and DMPlex for overlapping subdomains \cite{Knepley2009}. In our implementation, $A_{\bm{u}}$ is approximated by using the two-level additive Schwarz preconditioner from PETSc with PCHPDDM \cite{Jolivet2021} and SLEPc for GenEO.
\subsection{Dependency on the subdomain number $N$}
Scalability of GMRES-block triangular preconditioner is tested increasing the number of subdomains $N$. The subdomain size is set with $H/h=8$, and the overlapping width is set with $\delta/h=1$. Two cases are considered, 1) compressible and strongly permeable, 2) almost incompressible and weakly permeable. In the compressible and weakly permeable case, there comes a moderate scalable trend through $N=49$ and $N=64$ (Table 1, Left). 

\begin{table}[hbtp]
\begin{center}
  \small
  \begin{tabular}{c|c|c}
    \hline
&\multicolumn{1}{c|}{$\nu=0.3$} & \multicolumn{1}{|c}{$\nu=0.4999$}\\
&\multicolumn{1}{c|}{$\kappa=10^{-2}$} & \multicolumn{1}{|c}{$\kappa = 10^{-9}$}\\

         {$N$}         &iteration   & iteration    \\
    \hline
        $4$ &4   &  6     \\
    $9$ &5  & 9      \\
    $16$ &7  & 12     \\
    $25$ &7  & 14    \\
    $36$ &9  & 17      \\
    $49$ &11  & 19     \\
    $64$ &11  & 20    \\
    \hline
  \end{tabular}
\quad
  \begin{tabular}{c|c|c}
    \hline
&\multicolumn{1}{c|}{$\nu=0.3$} &\multicolumn{1}{c}{$\nu=0.4999$}\\
&\multicolumn{1}{c|}{$\kappa=10^{-2}$} &\multicolumn{1}{c}{$\kappa = 10^{-9}$}\\
         {$H/\delta$}         &iteration  & iteration  \\
             \hline
    8          & 52  & 56      \\
    $16$ & 59    & 52    \\
    $32$ & 60    & 48      \\
    $64$ & 78    &  50    \\
    \hline
  \end{tabular}
\caption{\textbf{Left}: Scalability of GMRES-block triangular preconditioner. Iteration counts for increasing number of subdomains $N$. 
Fixed $H/h=8, ~ \delta/h=1$, deflation $N$ = 15, $r_{\rm tol}=10^{-8}$. \textbf{Right}: $H/\delta$-dependency of GMRES-block triangular preconditioner. Iteration counts for increased overlapping. 
Fixed $N=16$ and $H/h=64$, $r_{\rm tol}=10^{-8}$.}
\end{center}
\end{table}
%

\subsection{Dependency on $H/\delta$}
$H/\delta$ dependency of GMRES-block triangular preconditioner is tested with the overlap changed from 1 to 8 with the domain size in each dimension $H/h=64$ and the number of subdomains $N=16$. The increment of overlap does reduce the GMRES iterations for the compressible and strongly permeable case. However, for the almost incompressible and weakly permeable case, the GMRES iterations increase with the change of $H/\delta$ from 32 to 16 and 8 (Table 1, Right).

\subsection{Dependency on permeability $k$} 
The robustness to $\kappa$ of the GMRES-block triangular preconditioner is tested with $\kappa$ changed from 1 to $10^{-9}$. Whether compressible or almost incompressible, GMRES iterations show biphasic pattern of decrease from $\kappa$=1 to $10^{-3}$ and increase when $\kappa$ changes towards ${\kappa}=10^{-9}$ (Table 2, Upper). Overall, the numerical scheme of GMRES-block triangular preconditioner show evident robustness in broad ranges of permeability both in compressible and almost incompressible regimes.

\subsection{Dependency on heterogeneous material properties of $\nu$ and $\kappa$}
To test the block preconditioner and GenEO two-level OAS solver for the displacement, the primary parameters of Poisson ratio $\nu$ and permeability $\kappa$ are treated as nonuniform. Figure 1(a) is in the pattern of checkboard and jump across subdomains and Figure 1(b) is in the pattern of jump along subdomains. 
Compressible and strongly permeable poroelasticity ($\nu$=0.3 and $\kappa=10^{-2}$)  is prescribed to yellow regions, and almost incompressible and weakly permeable poroelasticity ($\nu=0.4999$ and $\kappa=10^{-9}$) is prescribed to black regions. The proposed GMRES-block triangular preconditioner solver shows robustness to material heterogeneity with finite iterations for both patterns of non-uniformity (Table 2, Lower).
\begin{figure}[h]
\centering
$\begin{array}{c@{\hspace{0.05in}}c@{\hspace{0.05in}}c}
\includegraphics[angle=0, width=0.42\textwidth]{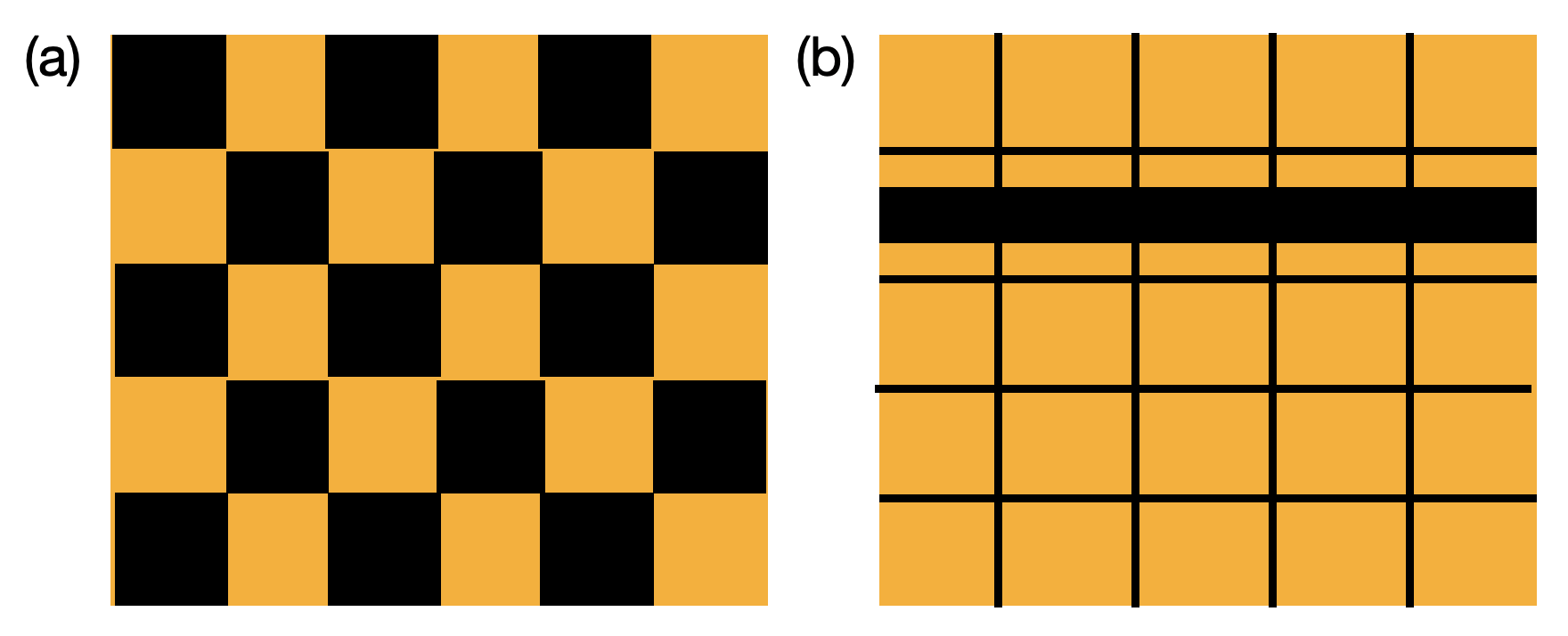}
\end{array}$
\qquad
  \begin{tabular}{c|c|c}
    \hline
&\multicolumn{1}{c|}{$\nu=0.3$} &\multicolumn{1}{c}{$\nu=0.4999$}\\
         {$k$}         &iteration     & iteration    \\
             \hline
    1          & 51  &  43    \\
    $10^{-1}$ & 19     & 17     \\
    $10^{-3}$ & 7    & 9     \\
    $10^{-5}$ & 14   & 9    \\
    $10^{-7}$ & 16   & 14    \\
    $10^{-9}$ & 16   & 47     \\
                 \hline
    &\multicolumn{2}{c}{$(\nu=0.3$ and $\kappa=10^{-2})$ vs.} \\
    &\multicolumn{2}{c}{$(\nu=0.4999$ and $\kappa=10^{-9})$}\\
                 \hline
    $\rm heter.-subdomain $ & 175      \\
    $\rm heter.-along $  &87        \\
    \hline
  \end{tabular}
\caption{(a) Jump across subdomains. (b) Jump along subdomains. \textbf{Table 2}: Robustness to $k$ of GMRES-block triangular preconditioner. Iteration counts for decreasing permeability. 
Fixed $N=16$ and $H/h=8$, deflation $N$ = 15, $r_{\rm tol}=10^{-10}$.}
\end{figure}
%
%

\section{Conclusion}

We proposed domain decomposition preconditioners for the saddle point problem of the three-field Biot model and performed numerical experiments for scalability and robustness in parameters. GMRES with block triangular preconditioner with two-level OAS with coarse space by GenEO for $\bm{u}$ and LU for $\bm{z}$ and $p$ is  scalable and robust in broad ranges of parameters ($\nu$, $\kappa$) and their heterogeneity. Future works are 1) some theoretical analysis on condition number bounds of proposed preconditioned systems, e.g. field-of-value analysis, 2) domain decomposition preconditioners for 3D poroelastic large deformation.  

\begin{acknowledgement}
The author gives thanks to Dr. Mingchao Cai for the introduction of Biot models and invaluable discussions. The author was partly supported by NSF DMS-1831950, and the virtual attendance to DD27 conference by Penn State University NSF travel fund.
\end{acknowledgement}

\end{document}